\title{
Two-spring 
topologies with minimal fabrication cost and required weighted force–resistance performance}
\author{Egor Makarenkov\footnote{Allen High School}}
\documentclass{article}
\usepackage{graphicx,amsmath}

\begin{document}
\newpage
\maketitle
\begin{abstract}
    Starting from a problem in elastoplasticity, we consider an optimization problem $C(c_1,c_2)=c_1+c_2\to \min$ under constraints $F_R^k(c_1,c_2)=a\cdot F^k(c_1,c_2)+b\cdot R^k(c_1,c_2)\ge 1$ and $F^k(c_1,c_2)\ge 1$, where both $F^k$ and $R^k$ non-linear, $a,b$ are constants, and $i\in\{1,2\}$ is an index. For each $(a,b)$ we determine which of the two values of $i\in\{1,2\}$ leads to the smaller minimum of the optimization problem. This way we obtain an interesting curve bounding the region where $k=1$ outperforms   $k=2$.

\end{abstract}



\section{Introduction}
\begin{figure}[h]
\begin{center}
    \includegraphics[scale=0.62]{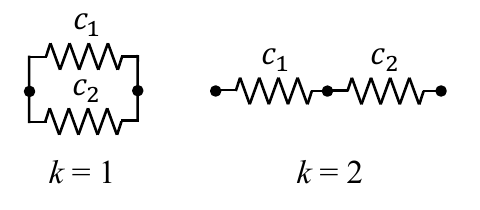}
\caption{All possible two-spring lattice configurations.
}
\label{configs}
\end{center}
\end{figure}

We consider the simplest network of two elastoplastic current conducting springs whose multi-functional performance $F_R^k$ is defined as a linear combination 
\begin{equation}\label{FR}
F_R^k(c_1,c_2)=a\cdot F^k(c_1,c_2)+b\cdot R^k(c_1,c_2)
\end{equation}
of the maximal force $F^k(c_1,c_2)$ that the system can withstand before deforming plastically and  electric resistance $R^k(c_1,c_2)$ of the system. The indices $k=1$ and $k=2$ correspond to parallel and serial connections of springs respectively, see Fig.~\ref{configs}.

\noindent The maximal force that the system can withstand before deforming plastically comes from maximal forces $c_1$ and $c_2$ (also known as springs' {\it elastic limits}) that springs 1 and 2 can withstand before deforming plastically. Specifically,
\begin{equation}\label{F1F2}
  F^1(c_1,c_2)=c_1+c_2, \qquad F^2(c_1,c_2)=\left\{\begin{array}{l}
  c_1,\quad {\rm if}\ c_1\le c_2,\\
  c_2,\quad {\rm if}\ c_2\le c_1.\end{array}\right.
\end{equation}
The formula for $F^2$ comes by observing that the force in a serial connection of springs doesn't increase (in response to stretching) after one of the springs reaches its elastic limit $c_i$ first, at which point the force in the system equals $c_i$. The formula for $F^1$ comes by observing that the force in parallel connection of springs stops increasing only after each of the springs reached their elastic limits. We refer the reader to \cite{Egor} for a comprehensive discussion of the response force formulas. Assuming that electric resistance of spring $i$ is $R_i=1/c_i$ (see \cite{mater01}), the resistances of serial and parallel connections of springs are given by
$$
   R^1(c_1,c_2)=\frac{1}{c_1+c_2},\qquad R^2(c_1,c_2)=\frac{1}{c_1}+\frac{1}{c_2},
$$
and the fabrication cost of the network is $C(c_1,c_2)=c_1+c_2$ (\cite{Sakshi,mater01}), the present paper addresses the optimization problem
\begin{equation}\label{op}
      c_1\ge 0,\ c_2\ge 0,\ F^k_R(c_1,c_2)\ge 1,\ F^k(c_1,c_2)\ge 1,\ k\in\{1,2\},\ c_1+c_2\to\min,
\end{equation}
whose optimal parameters $(c_1,c_2,k)$ and the minimal cost $c_1+c_2$ depend on parameters $a$ and $b$ in (\ref{FR}). The addition of $F^k(c_1,c_2)\ge 1$ means that we don't want to compromise the strength of springs' system while optimizing its multi-functional performance. We are most interested to find out whether one of the two systems of Fig.~\ref{configs} is always optimal regardless of $a$ and $b$ or $k(a,b)$ is a non-constant function of $a$ and $b$.  Physically, this means whether the optimal topology of spring system depends on multi-functional requirement or there is a universal topology which is optimal all the time.

\section{Reduction of the optimization problem to one parameter}\label{Sec2}

Observe that, when $k=1$ (the case of parallel springs), optimization problem (\ref{op}) can be rewritten through just one variable $x=c_1+c_2$, i.e.
\begin{equation}\label{op1reduced}
  x\ge 0, \ ax+\dfrac{b}{x}\ge 1, \ x\ge 1, \ x\to\min.
\end{equation}
A solution $x_*$ of optimization problem (\ref{op1reduced}) provides the minimal cost
\begin{equation}\label{cost2}
  C(c_{1*},c_{2*})=x_* 
\end{equation}
of optimization problem (\ref{op})
for any $c_{1*}\ge 0$ and $c_{2*}\ge 0 $ such that $x_*=c_{1*}+c_{2*}.$

\vskip0.2cm

\begin{figure}[h]
\begin{center}
    \includegraphics[scale=0.5]{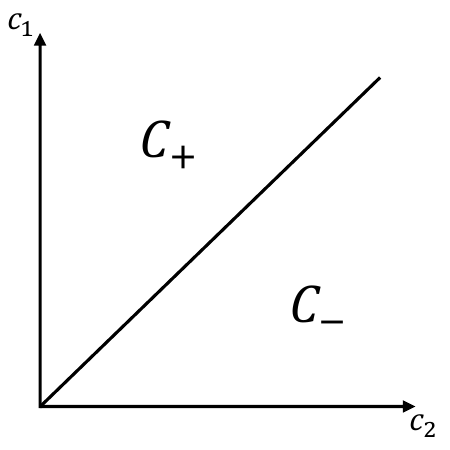}
\caption{Partition of the parameter space into ordering regimes that resolve the piecewise definition of the response force and enable dimension reduction.}
\label{figC+C-}
\end{center}
\end{figure}

\noindent It turns out that reduction to one variable is also possible in the case of $k=2$. First of all observe that 
\begin{itemize}
    \item[(i)] if we can find $\tilde F^2(c_1,c_2)$  that is always bigger (or equal) than $F^2(c_1,c_2)$, 
    \item[(ii)] if we can find 
$\tilde F_R^2(c_1,c_2)$ that is always bigger (or equal) than $F_R^2(c_1,c_2)$, 
    \item[(iii)] if we can find $\tilde C(c_1,c_2)$ that is always smaller (or equal) than $C(c_1,c_2)$,
    \item[(iv)] if an optimal value $(c_{1*},c_{2*})$ of the optimization problem
\begin{equation}\label{opnew}
      c_1\ge 0,\ c_2\ge 0,\ \tilde F^2_R(c_1,c_2)\ge 1,\ \tilde F^2(c_1,c_2)\ge 1,\ \tilde C(c_1,c_2)\to\min,
\end{equation}
satisfies 
\begin{equation}\label{CC}
C(c_{1*},c_{2*})=\tilde C(c_{1*},c_{2*}),
\end{equation}
\end{itemize}
then 
\begin{itemize}
    \item $(c_{1*},c_{2*})$ satisfies the domain of optimization problem (\ref{op}),
    \item an optimal solution of (\ref{op}) cannot be smaller than  $(c_{1*},c_{2*})$  because
    $$
      \tilde C(c_{1*},c_{2*})\le \tilde C(c_1,c_2)\le C(c_1,c_2),\quad \mbox{for all}\ c_1,c_2\ge 0. 
    $$
\end{itemize}
Therefore, if conditions (i)-(iv) hold, then $(c_{1*},c_{2*})$ is the optimal solution of problem (\ref{op}) (see  \cite[Proposition~1]{Egor} for an equivalent formulation of this observation). Though it is unlikely we can find $\tilde F^2$, $\tilde F_R^2$, and $\tilde C$ to satisfy condition (i)-(iv) in the entire domain $c_1\ge 0,$ $c_2\ge 0$, it appears we can do that separately in the domains
$$
C_+=\{(c_1,c_2):c_1\le c_2\}\quad \mbox{and}\quad C_-=\{(c_1,c_2):c_2\le c_1\}
$$
given by (\ref{F1F2}), see Fig.~\ref{figC+C-}. 

Indeed, since
$$
\renewcommand{\arraystretch}{2.5}
\begin{array}{l}
  R^2(c_1,c_2)=\dfrac{1}{c_1}+\dfrac{1}{c_2}\le \dfrac{2}{c_1},\quad C(c_1,c_2)=c_1+c_2\ge 2c_1\quad {\rm if}\ (c_1,c_2)\in C_+,\\
  R^2(c_1,c_2)=\dfrac{1}{c_1}+\dfrac{1}{c_2}\le \dfrac{2}{c_2},\quad C(c_1,c_2)=c_1+c_2\ge 2c_2 \quad {\rm if}\ (c_1,c_2)\in C_-
\end{array}
$$
we can take
\begin{eqnarray} 
  \tilde F^2(c_1,c_2)=c_1, \ \tilde F_R^2(c_1,c_2)=a c_1+b\dfrac{2}{c_1},\  \tilde C(c_1,c_2)=2c_1,\ \ {\rm if}\ (c_1,c_2)\in C_+,\label{op1}\\
    \tilde F^2(c_1,c_2)=c_2, \ \tilde F_R^2(c_1,c_2)=a c_2+b \dfrac{2}{c_2},\  \tilde C(c_1,c_2)=2c_2,\ \  {\rm if}\ (c_1,c_2)\in C_-,\label{op2}
\end{eqnarray}
to satisfy conditions (i)-(iii). Observe that condition (iv) will hold if solution $(c_{1*},c_{2*})$ of optimization problem (\ref{opnew}) in domains $C_+$ and $C_-$ belongs to the boundary of the corresponding domain, i.e. if 
\begin{equation}\label{c1*c2*}
c_{1*}=c_{2*}
\end{equation}
holds. Indeed, if (\ref{c1*c2*}) holds, then
$$
  \tilde C(c_{1*},c_{2*})=c_{1*}+c_{2*}=C(c_{1*},c_{2*})
$$
for each of the $\tilde C$ given by (\ref{op1}) and (\ref{op2}). But  
according to formulas (\ref{op1}) and (\ref{op2}) optimization problem (\ref{opnew}) reads as
\begin{equation}\label{op2reduced}
  x\ge 0, \ ax+\dfrac{2b}{x}\ge 1, \ x\ge 1, \ 2x\to\min
\end{equation}
in each of the domains $C_+$ and $C_-$, if $x$ is chosen to denote $c_1$ in the domain $C_+$ and $c_2$ in the domain $C_-$. Therefore, to solve optimization problem (\ref{op}) for $k=2$, we have to find the optimal value $x_*$ of (\ref{op2reduced}) and take
$$
  c_{1*}=x_*,\ \ c_{2*}=x_*, 
$$
obtaining the cost
\begin{equation}\label{cost1}
    C(c_{1*},c_{2*})=2x_*.
\end{equation}

\noindent In the next section of the paper we solve the reduced optimization problems (\ref{op1reduced}) and (\ref{op2reduced}) for each $a\ge 0$ and $b\ge 0$ (obtaining optimal values of $x_*^{k=1}$ and $x_*^{k=2}$) and conclude that

\vskip0.2cm

\centerline{$k=1$ wins if $x_*^{k=1}>2 x_*^{k=2}$ and $k=2$ wins if 
$2 x_*^{k=2}>x_*^{k=1},$} 

\vskip0.2cm

\noindent which inequalities come from (\ref{cost2}) and (\ref{cost1}).

\section{Solutions of 
reduced optimization problems}

\begin{figure}[h]
\begin{center}
    \includegraphics[scale=0.57]{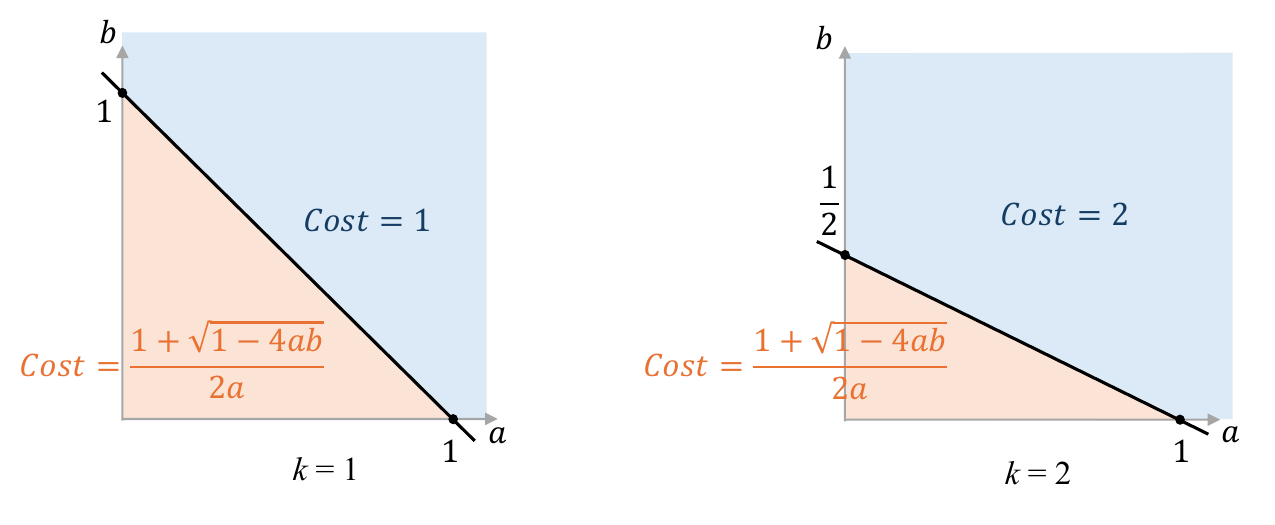}
\caption{Visualization of the reduced minimal cost expression (\ref{eqcost}). 
}
\label{figeqcost}
\end{center}
\end{figure}

For fixed $a\ge 0$ and $b\ge 0$, the solution of  optimization problems (\ref{op1reduced}) and (\ref{op2reduced}) is the minimal solution of the system of inequalities
\begin{eqnarray}
ax+\dfrac{kb}{x}&\ge& 1,\label{inq1}\\
x&\ge& 1,\label{inq2}
\end{eqnarray}
 where $k=1$ or $k=2$ accordingly.
Solution of (\ref{inq1}) is 
$$
\begin{array}{ll}
  \left(-\infty,x_1\right]\cup \left[x_2,\infty\right), & {\rm if}\ 1-4kab\ge 0,\\
  (-\infty,\infty),& {\rm if}\  1-4kab \le 0,
\end{array}
$$
where 
$$
x_1=\dfrac{1-\sqrt{1-4kab}}{2a},\quad x_2=\dfrac{1+\sqrt{1-4kab}}{2a}.
$$
Therefore, when $1-4kab\le 0$, the minimal solution of (\ref{inq1})-(\ref{inq2}) is $x=1$. When $1-4kab\ge 0$, the minimal solution of (\ref{inq1})-(\ref{inq2}) is
$$
\left\{
\begin{array}{ll}
   1,& {\rm if}\  1\in \left(-\infty,x_1\right]\ {\rm or}\ 1\in \left[x_2,\infty\right),\\
   x_2,& {\rm if}\ 1\in \left(x_1,x_2\right).
\end{array}\right.
$$
In summary, the minimal solution of (\ref{inq1})-(\ref{inq2}) is 
\begin{equation}
\left\{
\begin{array}{ll}
   x_2,& {\rm if}\ 1\in \left(x_1,x_2\right)\ {\rm and}\ 1-4kab\ge 0,\\
      1,& {\rm otherwise}.
\end{array}\right.
\end{equation}
Observing that, when $a\ge 0$ and $b\ge 0$, the system of inequalities 
$1\in \left(x_1,x_2\right)$ and $1-4kab\ge 0$ is equivalent to just $a+kb-1<0$,
we conclude that the minimal cost of optimization problems (\ref{op1reduced}) and (\ref{op2reduced}) is
\begin{equation}\label{eqcost}
\left\{
\begin{array}{ll}
   k\dfrac{1+\sqrt{1-4kab}}{2a}, & {\rm if}\ a+kb-1<0,\\
      k, & {\rm otherwise},
\end{array}\right.
\end{equation}
with $k=1$ and $k=2$ accordingly. Formula (\ref{eqcost}) is illustrated in Fig.~\ref{figeqcost}.

\subsection{Determining the winning configurations of springs}

\begin{figure}[h]
\begin{center}
    \includegraphics[scale=0.6]{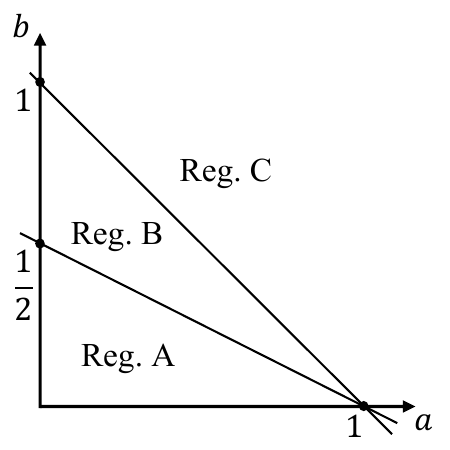}
\caption{Subdivision of the weight space $(a,b)$ indicating which constraint is active and which reduced cost formula governs the optimization.}
\label{figABC}
\end{center}
\end{figure}

\noindent Based on formula (\ref{eqcost}) for the cost, the space of parameters $(a,b)$ splits into 3 regions separated by the line $a+kb-1=0$ with $k=1$ and $k=2$:
\begin{eqnarray*}
    \mbox{region A:} && \{(a,b):\ a\ge 0,\ b\ge 0,\ a+2b-1<0\},\\
    \mbox{region B:} && \{(a,b):\ a\ge 0,\ b\ge 0,\ a+2b-1\ge 0\ {\rm and}\ a+b-1<0\},\\
    \mbox{region C:} && \{(a,b):\ a\ge 0,\ b\ge 0,\ a+b-1\ge 0\},
\end{eqnarray*}
see Fig.~\ref{figABC}.

In what follows we go over each of the regions and determine the winning configuration of springs in each of the regions.

\vskip0.2cm

\noindent {\bf Region A.} According to the summary in Fig.~\ref{figeqcost}, $k=1$ outperforms $k=2$ (i.e. the configuration $k=1$ delivers smaller minimal cost compared to configuration $k=2$) in Region~A, if 
$$
  \dfrac{1+\sqrt{1-4ab}}{2a}<2\dfrac{1+\sqrt{1-8ab}}{2a},
$$
which does hold for all $(a,b)$ from region $A.$

\vskip0.2cm

\begin{figure}[h]
\begin{center}
    \includegraphics[scale=0.6]{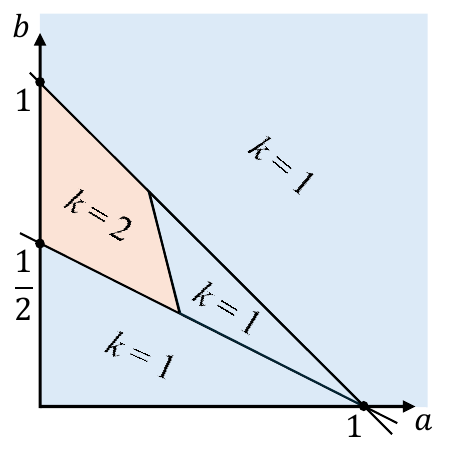}
\caption{Phase diagram that, for each $(a,b)$, identifies the topology (parallel $k=1$ or serial $k=2$) that solves optimization problem (\ref{op}), i.e. achieves the lowest fabrication cost.}
\label{figsum}
\end{center}
\end{figure}

\noindent {\bf Region B.}  The case $k=1$ outperforms $k=2$ in Region~B, if
$$
  \dfrac{1+\sqrt{1-4ab}}{2a}<2,
$$
which does hold in a part of region~B, that is highlighted in Fig~\ref{figsum}.

\vskip0.2cm

\noindent {\bf Region C.}  The case $k=1$ outperforms $k=2$ in Region~C, because $1<2$.

\vskip0.2cm

\noindent Our findings can be summarized as the following statement.

\vskip0.2cm

\noindent {\bf Proposition.} For each values of $a\ge 0$ and $b\ge 0$, the minimal cost in optimization problem (\ref{op}) is attained at $k=1$ except for the special region
\begin{eqnarray*}
    \mbox{region B2:} && \left\{(a,b):\ a+2b-1\ge 0,\ a+b-1<0,\   \dfrac{1+\sqrt{1-4ab}}{2a}>2\right\},
\end{eqnarray*}
where the minimal cost of (\ref{op}) is attained at $k=2.$

\section*{Acknowledgments}
\noindent The author is grateful to his mentors, Sakshi Malhotra (Department of Mathematics and Statistics, Montgomery College, Maryland) and Yang Jiao (School for Engineering of Matter, Transport and Energy, Arizona State University), for guiding him throughout the project, and to his Advanced Pre-Calculus teacher, Mrs. Williams, for building a strong foundation in the analysis of polynomial and rational functions, that supported this work.

\bibliographystyle{plain}

\noindent {\small Egor Makarenkov, Allen High School, 300 Rivercrest Boulevard, Allen, TX 75002, egordaporg@gmail.com, egor.makarenkov@student.allenisd.org}

\vskip0.2cm

\noindent Egor Makarenkov is a junior at Allen High School.

\end{document}